\documentclass[11pt,twosides]{amsart}

\usepackage[dvips]{color}
\usepackage{amsfonts}
\usepackage{amssymb, latexsym, amsmath, pb-diagram}
\usepackage{pstricks-add}
\usepackage{lamsarrow, pb-lams}
\usepackage{multicol}
\usepackage{bm}
\usepackage[mathscr]{eucal}
\usepackage{xy}
\usepackage{epic,eepic}
\xyoption{all}

\def\qb{\hfill $\Box$}

\begin{document}
\title[Torsion]{A moment-angle manifold whose cohomology has torsion}
\author[X.Li \& G.Wang]{Xiaomeng Li and Gefei Wang}
\date{}
\address{School of Mathematical Science, Nankai University, Tianjin 300071, P.~R.~China}
\email{1412640@mail.nankai.edu.cn}
\email{lxmnku1992@163.com}
\thanks{This project is supported by NSFC  No.11471167 and No.11761072}
\keywords{moment-angle manifold, stellar subdivision, full subcomplex, missing face and simplicial complement}
\subjclass[2010]{13F55, 05A19, 05E40, 52B05, 52B10}
\maketitle
\begin{abstract}
  In this paper we give a method to construct moment-angle manifolds whose cohomology
has torsion. We also give method to describe the corresponding simplicial sphere
by its non-faces.
\end{abstract}

\section{Introduction}
  Corresponding to every abstract simplicial complex $L$ on the vertex set $[m]=\{1,2,\cdots, m\}$,
there are the real and complex moment-angle complexes $\mathbb{R}\mathcal{Z}_{L}$ and $\mathcal{Z}_{L}$
({\it cf.} \cite{BP02, BP15}).
They are defined as
\begin{align*}
\mathbb{R}\mathcal{Z}_{L} = & \bigcup_{\sigma \in L}\prod_{i \in \sigma}D_i^1 \times \prod_{i \notin \sigma}S_i^0
  \subseteq D_1^1\times D_2^1\times \cdots \times D_m^1 \\
\mathcal{Z}_{L} = & \bigcup_{\sigma \in L}\prod_{i \in \sigma}D_i^2 \times \prod_{i \notin \sigma}S_i^1
  \subseteq D_1^2\times D_2^2 \times \cdots \times D_m^2.
\end{align*}

The cohomology groups of $\mathbb{R}\mathcal{Z}_{L}$ and $\mathcal{Z}_{L}$ are given by
Hochster's theorem:

{\bf Theorem \cite{B02, BP02, BP15}}
  {\it Let $L$ be a simpicial complex on the vertex set [m], then
  \begin{align*}
  H^{*}(\mathbb{R}\mathcal{Z}_L) &\cong \bigoplus_{I \subset [m]}\widetilde{H}^{*-1}(L|_I) \\
  H^{*}(\mathcal{Z}_L) &\cong \bigoplus_{I \subset [m]}\widetilde{H}^{*-|I|-1}(L|_I)
  \end{align*}
where $L|_I$ is the full subcomplex of $L$ on subset $I$  and $I$ runs over all the subsets of $[m]$.
}

  From \cite{BPR07, BR01, Cai15} we know that both $\mathbb{R}\mathcal{Z}_{L}$ and $\mathcal{Z}_{L}$
are topological manifolds if $L$ is a simplicial sphere, referred to as moment-angle manifolds.
Furthermore if $L$ is a polytopal sphere (the boundary complex of a simplicial polytope),
then $\mathcal{Z}_{L}$ is a transverse intersection of real quadratic hypersurfaces
({\it cf. \cite{BM06}}), while both $\mathbb{R}\mathcal{Z}_L$ and $\mathcal{Z}_L$ are framed
differentiable manifolds.

F.Bosio and L.Meersseman in \cite{BM06} announced that the cohomology groups of
differentiable complex moment-angle manifolds may have any torsion $\mathbb{Z}/m$.
Furthermore if $L$ is $\mathbb{Z}/2$ colourable, L.Cai, S.Choi and H. Park in \cite{CC16, CP13} proved that
the small cover under $\mathbb{R}\mathcal{Z}_L$ may have any torsion $\mathbb{Z}/m$.

  From  Hochster's theorem, it is easy to construct a moment-angle complex whose cohomology has torsion.
But it is harder to construct such moment-angle manifolds
(at least, the cohomology of all the moment-angle manifolds corresponding to dimensional 1, 2 and 3
simplicial spheres are torsion free ({\it cf.} \cite{BM06} Corollary 11.1).

  Based on  Hochster's theorem, our goal is to find a simplicial complex $K$
whose cohomology has torsion and $K$ is embedded
in a polytopal sphere $L$ as a full subcomplex.
Then both the real and complex moment-angle complexes corresponding to $L$ are differentiable manifolds
and the cohomology of $\mathbb{R}\mathcal{Z}_L$ and $\mathcal{Z}_L$ have $\widetilde{H}^{*}(K)$
as a summand and then have torsion.

{\bf Theorem 3.2 \space (Construction)}
 {\it Let $K$ be a subcomplex (not a full subcomplex) of a simplicial sphere $L_0$ on the vertex set $[m]$,
$\mathbb{M} = \{\sigma_1,\sigma_2,\cdots,\sigma_{s}\}$ be the set of missing faces of $K$, which are also simplices of $L_0$.
As following, make stellar subdivisions at $\sigma_1,\sigma_2,\cdots,\sigma_{s}$ on $L_0$ one by one
\[\begin{array}{cccccccccccc}
L_1 =ss_{\sigma_1}L_0,  & L_2 =ss_{\sigma_2}L_1, & \cdots,&  L_{s} =ss_{\sigma_{s}}L_{s-1}.
\end{array}\]
Then $K$ becomes a full subcomplex of $L_s$, $K = L_s|_{[m]}$.}

In fact after making stellar subdivision on a polytopal (simplicial) sphere, it
is still polytopal (simplicial) (see \cite{ES74}). If $L_0$ is also a polytopal sphere, we thus obtain a polytopal sphere
$L_s$  by Theorem $3.2$ such that $K$ is a full subcomplex of $L_s$. The real and complex
moment-angle complexes corresponding to $L_s$ are differentable manifolds. By Hochster's theorem
both $H^*(\mathbb{R}\mathcal{Z}_{L_s})$ and $H^*(\mathcal{Z}_{L_s})$ have torsion if $\widetilde{H}^*(K)$ has torsion.

  At last in section 4, we give a differentiable moment-angle manifold whose cohomology
has $\mathbb{Z}/3$ as a summand. This  is done as follows:

  Triangulate the $mod$ $3$ Moore space $K$  which has $8$ vertices, $17$ $2$-dimensional facets and $22$ missing faces
(see figure 3). It can be embedded in $\partial\Delta^7=L_0$.
After making $22$ stellar subdivisions on it,  $K$ becomes a full subcomplex
of the polytopal sphere $L_{22}$.
Then $L_{22}$ is a $6$ dimensional polytopal sphere with $30$ vertices. $\mathcal{Z}_{L_{22}}$ is a
$37$ dimensional differentiable manifold and $H^{11}(\mathcal{Z}_{L_{22}})$ has
$\widetilde{H}^2(K)=\mathbb{Z} / 3$ as a summand.

It is notable that F.Bosio and L.Meersseman's construction in \cite{BM06} (Theorem $11.12$)
applied to the same example does not give a moment-angle manifold whose cohomology has $\mathbb{Z}/3$ torsion.

The authors are grateful to professor Zhi. L$\ddot{u}$ for his helpful suggestion during this
research.
This work was done under the supervision of professor Xiangjun. Wang.

\section{Simplicial Complement}

  An abstract simplicial complex $K$ on the vertex set $I$ is a collection of simplices that satisfies:
for any simplex (face) $\sigma\in K$, all of its proper subsets (proper faces) are simplices of $K$.

  An abstract simplicial complex $K$ could also be given by all of its {\it non-faces}
\[
\mathbb{A}=2^{I}\setminus K
\]
and
\[
K=2^I \setminus \mathbb{A}
\]
that satisfies: if $\sigma\in\mathbb{A}$ is  not a simplex of $K$ and $\sigma'\supset\sigma$
then $\sigma'\in\mathbb{A}$ is not a simplex of $K$.

  A simplex $\sigma=(i_1,i_2,\cdots,i_k)\in 2^{I}$ is called a {\it missing face} (or minimal non-face) of $K$ if it is not
a face of $K$, but all of its proper subsets are faces of $K$, i.e. $\sigma \notin K$ but every $\sigma_j =(i_1,\cdots,\widehat{i_j},\cdots,i_k) \in K$ , $j=1,2,\cdots, k$.
An abstract simplicial complex could also be given by its set of missing faces
\[
\mathbb{M}=\{\sigma\in 2^{I}|\mbox{$\sigma$ is a missing face of $K$}\}
\]
and
\[
K=\{\tau \in 2^{I}| \mbox {$\tau$ does not contain any $\sigma \in \mathbb{M}$ }\}.
\]
A subset $\sigma'$ of $I$ is not a simplex of $K$ if and only if it contains a missing
face $\sigma\in\mathbb{M}$ as a subset.

{\bf Definition 2.1}
{\it Let $K$ be a simplicial complex on the vertex set $I$ and $\mathbb{M}$, $\mathbb{A}$ be the sets of
missing faces and non-faces of $K$ respectively. We define a simplicial complement of $K$, denoted by
\[
\mathbb{P}=\{\sigma_1, \sigma_2, \cdots, \sigma_s\},
\]
to be a collection of non-faces that includes all the missing faces $\mathbb{M}$ i.e.
\[
\mathbb{M} \subseteq \mathbb{P} \subseteq \mathbb{A}.
\]}

Similar to the set of missing faces $\mathbb{M}$, given a simplicial complement $\mathbb{P}$ (collection of non-faces) on the vertex set $I$, one can
obtain a simplicial complex $K_{\mathbb{P}}$ on $I$ by:
 \[
 K_{\mathbb{P}}(I)=\{\tau \subset I|\mbox{$\tau$ does not contain any $\sigma_i \in \mathbb{P}$} \} \tag{2.2}
 \]
 or by all of its non-faces
 \[
2^{I} \setminus K_{\mathbb{P}}(I)=\{\tau \subset I|\mbox{$\tau$ contains a $\sigma_i \in \mathbb{P}$} \}\tag{2.3}.
 \]

   A subset $\sigma$ of $I$ is
not a simplex of $K_{\mathbb{P}}(I)$ if and only if it contains a non-face $\sigma_{i}$ in the
simplicial complement $\mathbb{P}$.

  {\bf Definition 2.4}
  {\it Let $\mathbb{P}, \mathbb{P}'$ be two simplicial complements on the vertex set $I$, if they can obtain the same
simplicial complex i.e. $K_{\mathbb{P}}(I)=K_{\mathbb{P}'}(I)$,
we say that $\mathbb{P}$ and $\mathbb{P}'$ are equivalent, denoted by $\mathbb{P}\simeq \mathbb{P}'$.}

It is easy to see that:
 Two simplicial complements $\mathbb{P}$\textcolor{red}{,} $\mathbb{P}^{'}$ on $I$ are equivalent
if and only if for every non-face $\sigma \in \mathbb{P}$ there exists a $\sigma' \in \mathbb{P}'$ such that
$\sigma' \subseteq \sigma$ and for every non-face $\sigma' \in \mathbb{P}'$ there exists a
$\sigma \in \mathbb{P}$ such that $\sigma \subseteq \sigma'$.

{\bf Proposition 2.5}
 {\it Let $\mathbb{P}=\{\sigma_1,\sigma_2,\cdots,\sigma_s\}$ be a simplicial complement of $K$ on $I$. For a
non-face $\sigma_j\in\mathbb{P}$ if there exists a $\sigma_i\in \mathbb{P},i \neq j$ such that
$\sigma_i\subseteq \sigma_j$, then we can remove $\sigma_j$  from $\mathbb{P}$ and the resulting simplicial complement
\[
\mathbb{P}'=\{\sigma_1,\sigma_2,\cdots,\widehat{\sigma}_j,\cdots,\sigma_s\}
\]
is equivalent to $\mathbb{P}$. In this case we call that $\mathbb{P}$ is reduced to $\mathbb{P}'$.}

Every simplicial complement of $K$ could be reduced to the set of missing faces by removing all the
larger non-faces.

{\bf Example 1} The simplicial complex $K$ is determined by the maximal simplices $(1,3)$, $(2,3)$, $(1,2,4)$, $(1,2,5)$, $(1,4,5)$, $(2,4,5)$ and their
proper subsets on the vertex set $[5]=\{1,2,3,4,5\}$ (see Figure 1)

\begin{center}
\setlength{\unitlength}{0.7mm}
\begin{picture}(100,45)(-50,-15)
\put(10,25){\line(0,-1){40}}
\put(30,5){\line(-1,1){20.2}}
\put(30,5){\line(-1,-1){20.2}}
\put(-30,-5){\line(81,61){40}}
\put(-30,-5){\line(83,-21){40}}
\put(-30,-5){\line(6,1){25}}
\put(10,25){\line(-15,-26){15}}
\put(10,-15){\line(-20,19){15}}
\put(10,27){$1$}
\put(12,-17){$2$}
\put(31,5){$3$}
\put(-33,-5){$4$}
\put(-7,-6){$5$}
\end{picture}
\end{center}

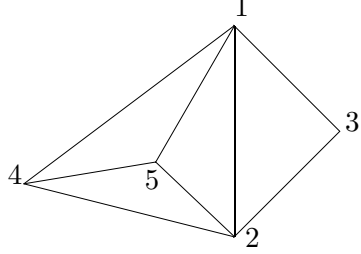
\centerline{Figure 1}

\begin{align*}
\mathbb{P}=
\left\{\begin{array}{lllll}
     \sigma_1=(1,2,4,5), \sigma_2=(1,2,3), \sigma_3=(3,4),\\
     \sigma_4=(3,5), \sigma_5=(1,3,4), \sigma_6=(3,4)
            \end{array}\right\}
\end{align*}
is a simplicial complement of $K$ on vertex set $[5]$ where $\sigma_3=(3,4)=\sigma_6$ appeared twice
and $\sigma_3=(3,4)\subset \sigma_5=(1,3,4)$. So $\sigma_5=(1,3, 4)$ and $\sigma_6=(3,4)$ could be removed from $\mathbb{P}$
to reduce to the set of missing faces $\mathbb{M}=\{(1,2,4,5), (1,2,3), (3,4), (3,5)\}$.

The readers should be aware that the  empty simplex
$\{\emptyset\}$ (only the empty set is a simplex) is different from the empty complex $\emptyset$
(the empty set is not a simplex of $\emptyset$).
$\mathbb{M}=\{(1), (2), \cdots, (m)\}$ is the set of missing faces of the empty simplex
$\{\emptyset\}$ while $\mathbb{M}_1=\{\emptyset\}$ is the set of missing faces of the empty
the empty complex $\emptyset$.

  Let $K$ be a simplicial complex on the vertex set $I$ and $\sigma$ be a simplex of $K$.
The $link$ and $star$ of $\sigma$ are defined to be the simplicial complexes
\begin{align*}
link_K\sigma = & \{\tau\in K| \sigma\cup\tau\in K, \sigma\cap\tau=\emptyset\},
  & star_K\sigma = & \{\tau\in K|\sigma\cup\tau\in K\}.
\end{align*}
The interior (open) $star$ is defined to be a subset of $2^{I}$ ({\it cf.} \cite{Mks84} \S 62 p.371)
\[
Intstar_K\sigma=\{\tau\in K|\sigma\subset \tau\}
\]
and the boundary of $star$ is the simplicial complex
\[
\partial star_K\sigma=star_K\sigma\setminus Intstar_K\sigma=\{\tau\in K|\sigma\cup\tau\in K, \sigma\not\subset \tau\}.
\]

Let $K_1$ and $K_2$ be two simplicial complexes on vertex set $I$ and $J$, where $I \cap J = \emptyset$.
The join of $K_1$ and $K_2$ is defined to be the simplicial complex on vertex set $I \cup J$
\[
K_1 \ast K_2=\{\sigma \cup \tau \in 2^{ I \cup J} | \sigma \in K_1 , \tau \in K_2\}.
\]

Let $\mathbb{P}=\{\sigma_1,\sigma_2,\cdots,\sigma_s\}$ be a simplicial complement of $K$ on the vertex set $I$ and $\sigma \in K$
be a simplex. We difine
\[
\mathbb{P}-\sigma = \{\sigma_1 \setminus \sigma,\sigma_2 \setminus \sigma,\cdots,\sigma_s \setminus \sigma\}
\]
which is a sequence of subsets on $I \setminus \sigma$.

{\bf Lemma 2.6} {\it Let $\mathbb{P}=\{\sigma_1,\sigma_2,\cdots,\sigma_{s}\}$
be a simplicial complement of $K$ on the vertex set $I$. Then
 \begin{enumerate}
\item
$
\mathbb{P}-\sigma = \{\sigma_1 \setminus \sigma,\sigma_2 \setminus \sigma,\cdots,\sigma_s \setminus \sigma\}
$
is a simplicial complement of $link_{K}\sigma$ on the vertex set $I\setminus \sigma$, i.e. by $(2.2)$
\[
link_K\sigma =K_{\mathbb{P}-\sigma}(I\setminus \sigma)=\{\tau\subset (I\setminus\sigma)|
  \mbox{$\sigma$ does not contain any $\sigma_i\setminus\sigma\in\mathbb{P}-\sigma$}\}.
\]
\item If we consider $\mathbb{P}-\sigma$ as a sequence of non-faces on the vertex set $I$, then it is a simplicial complement
of $star_K\sigma$ on I, i.e. by $(2.2)$
\[
star_{K}\sigma=K_{\mathbb{P}-\sigma}(I)=\{\tau\subset I|
  \mbox{$\sigma$ does not contain any $\sigma_i\setminus\sigma\in\mathbb{P}-\sigma$}\}.
\]
\end{enumerate}
}

  {\bf Proof:} We prove this Lemma by showing that they have the same non-faces
\begin{align*}
2^{I\setminus \sigma} \setminus link_K\sigma = & 2^{I\setminus \sigma} \setminus K_{\mathbb{P}-\sigma}(I\setminus \sigma)
  =  \{\tau\subset (I\setminus\sigma)|\mbox{$\tau$ contains a $\sigma_i\setminus\sigma\in\mathbb{P}-\sigma$}\}
\end{align*}
  and $$2^{I} \setminus star_K\sigma =2^{I} \setminus K_{\mathbb{P}-\sigma}(I)
=  \{\tau\subset I|\mbox{$\tau$ contains a $\sigma_i\setminus\sigma\in\mathbb{P}-\sigma$}\}. $$
  \begin{enumerate}
  \item From its definition, we know that
a simplex $\tau$ on the vertex set $I\setminus \sigma$ is not a simplex of $link_K\sigma$ if and only if
$\sigma \cup \tau$ is not a simplex of $K$. In other words, there exists a
$\sigma_i \in \mathbb{P}$ such that $\sigma_i \subseteq \tau \cup \sigma$.
This is equivalent to say that $\sigma_i \setminus \sigma \subseteq \tau \cup \sigma \setminus \sigma = \tau$,
every non-face $\tau \in 2^{I\setminus \sigma} \setminus link_K\sigma$ is a non-face of $K_{\mathbb{P}-\sigma}(I\setminus \sigma)$, i.e. $\tau \in 2^{I\setminus \sigma} \setminus K_{\mathbb{P}-\sigma}(I\setminus \sigma)$,
so $$2^{I\setminus \sigma} \setminus link_K\sigma \subseteq 2^{I\setminus \sigma} \setminus K_{\mathbb{P}-\sigma}(I\setminus \sigma).$$

  \item If a simplex $\tau$ on the vertex set $I \setminus \sigma$ contains a $\sigma_i \setminus \sigma$,
  then $\tau \cup \sigma \supseteq (\sigma_i \setminus \sigma) \cup \sigma \supseteq  \sigma_i$, so such $\tau$ is not a simplex of $link_K\sigma$.
  This is equivalent to say that every non-face $\tau \in 2^{I \setminus \sigma} \setminus K_{\mathbb{P}-\sigma}(I\setminus \sigma)$ is a non-face of $link_K\sigma$,
  i.e. $\tau \in 2^{I\setminus \sigma} \setminus link_K\sigma$, so
 \[
2^{I\setminus \sigma} \setminus K_{\mathbb{P}-\sigma}(I\setminus \sigma)\subseteq 2^{I\setminus \sigma} \setminus link_K\sigma.
\]
  \end{enumerate}

Thus $\mathbb{P}-\sigma = \{\sigma_1 \setminus \sigma,\sigma_2 \setminus \sigma,\cdots,\sigma_s \setminus \sigma\}$
is a simplicial complement of $link_{K}\sigma$ on the vertex set $I \setminus \sigma$.

  Similarly, if we consider $\mathbb{P}-\sigma$ as a simplicial complement on the vertex set $I$, then
\[
2^{I} \setminus star_K\sigma =2^{I} \setminus K_{\mathbb{P}-\sigma}(I)
=  \{\tau\subset I|\mbox{$\tau$ contains a $\sigma_i\setminus\sigma\in\mathbb{P}-\sigma$}\}.
\]
The Lemma follows. \qb

{\bf Example 2} In Example 1, the $link$ of the simplex $(1,2)$ is two vertices $link_K(1,2)=\{ (4), (5)\}$
and $star_K(1,2)$ is composed by two $2$-simplices $(1,2,4)$, $(1,2,5)$ and its proper subsets.
\begin{align*}
\mathbb{M}-(1,2)= & \{(1,2,4,5) \setminus(1,2),\ \ (1,2,3) \setminus (1,2),\ \ (3,4)\setminus(1,2),\ \ (3,5)\setminus(2)\}\\
=&\{(4,5),\ (3),\ (3,4),\ (3,5)\}\\
\simeq & \{(4,5), (3)\}
\end{align*}
is a simplicial complement of $link_K(1,2)$ on the vertex set $\{3,4,5\}$. Consider $\mathbb{M}-\sigma$
as a sequence of non-faces on the vertex set
$[5]=\{1,2,3,4,5\}$, it becomes the simplicial complement of $star_K(1,2)$.

 Let $\mathbb{P}_1=\{\sigma_{1}, \sigma_{2}, \cdots, \sigma_{s}\}$ and
$\mathbb{P}_2=\{\tau_{1},\tau_{2}, \cdots, \tau_{t}\}$ be the simplicial complements of $K_1$ and $K_2$
on the vertex set $I$.  We define their
join $\mathbb{P}_1\ast\mathbb{P}_2$ to be
\[
\mathbb{P}_1\ast\mathbb{P}_2=
\{\sigma_i\cup\tau_j|\sigma_i\in\mathbb{P}_1, \tau_j\in\mathbb{P}_2\},
\]
which is a sequence of subsets on $I$.

{\bf Lemma 2.7}  {\it Let $K_1$ and $K_2$ be two simplicial  complexes on the vertex set $I$,
$\mathbb{P}_1=\{\sigma_{1}, \sigma_{2}, \cdots, \sigma_{s}\}$ and
$\mathbb{P}_2=\{\tau_{1},\tau_{2}, \cdots, \tau_{t}\}$ be simplicial complements of
$K_1$ and $K_2$ respectively. Then
$\mathbb{P}_1 \ast \mathbb{P}_2=\{ \sigma_{i} \cup \tau_{j} | \sigma_{i} \in \mathbb{P}_1, \tau_{j} \in \mathbb{P}_2\}$
is a simplicial complement of $K_1 \cup K_2$ on the vertex set $I$,
\[
K_1\cup K_2=K_{\mathbb{P}_1\ast\mathbb{P}_2}(I)=\{\tau\subset I|\mbox{$\tau$ does not contain any $\sigma_i\cup\tau_j\in
  \mathbb{P}_1\ast\mathbb{P}_2$}\}.
\]}

{\bf Proof:} We prove this Lemma in the same way as the proof of Lemma $2.6$.
\begin{enumerate}
\item It is easy to see that a simplex $\tau$ on the vertex set $I$ is not a simplex of $K_1 \cup K_2$
if and only if it is not a simplex of either $K_1$ or $K_2$. This implies that there exists a
$\sigma_{i} \in \mathbb{P}_1$ such that $\sigma_i\subseteq \tau$  and also exists a
$ \tau_{j} \in \mathbb{P}_2$ such that $\tau_{j} \subseteq \tau$. This is equivalent to say that $\sigma_{i} \cup \tau_{j}\subseteq \tau$,
every non-face $\tau$ of $K_1 \cup K_2$ contains a $\sigma_i\cup\tau_j\in\mathbb{P}_1 \ast \mathbb{P}_2$, so
$$2^{I}\setminus  K_1 \cup K_2 \subseteq 2^{I}\setminus K_{\mathbb{P}_1 \ast \mathbb{P}_2}(I).$$

\item If a simplex $\tau$ on $I$ contains a non-face $\sigma_i\cup\tau_j\in\mathbb{P}_1 \ast \mathbb{P}_2$,
then $\sigma_i \subseteq \tau$ and $\tau_j \subseteq \tau$. This $\tau$ is neither a simplex of $K_1$ nor
a simplex of $K_2$, so
$$2^{I}\setminus K_{\mathbb{P}_1 \ast \mathbb{P}_2}(I) \subseteq 2^{I}\setminus K_1 \cup K_2.$$
\end{enumerate}
The Lemma follows.\qb

{\bf Corollary 2.8} {\it If the simplicial complement $\mathbb{P}$ is equivalent to $\mathbb{P}'$,
then for any simplex $\sigma$ and simplicial complement $\mathbb{P}_2$
\begin{align*}
\mathbb{P}-\sigma\simeq  & \mathbb{P}'-\sigma,
  & \mathbb{P}\ast\mathbb{P}_2\simeq & \mathbb{P}'\ast \mathbb{P}_2.
\end{align*}}
\qb

  Let $\sigma$ be a simplex of a simplicial complex $K$ on $[m]$. The stellar subdivision
at $\sigma$ on $K$ is defined to be the union of the simplicial complexes $K\setminus Intstar_K\sigma$
and the cone $cone\partial star_K\sigma$ along their boundary $\partial star_K\sigma$, denoted by
\[
ss_\sigma K =(K\setminus Intstar_K\sigma)\cup(cone\partial star_K\sigma),
\]
where
\[
cone\partial star_K \sigma= (m+1) \ast \partial star_K \sigma.
\]

  After stellar subdivision, one more vertex is added which is the vertex of the cone.({\it cf.} \cite{BM06})

In \cite{BP15} (Definition $2.7.1$), the stellar subdivision is defined to be
$$ss_\sigma K=(K\setminus star_K\sigma)\cup(cone\partial star_K\sigma),$$
where $K\setminus star_K\sigma$ is not a simplicial complex and $\sigma$ can be any simplex of $K$.
Note that
\[
K\setminus star_K\sigma= (K\setminus Intstar_K\sigma)\setminus\partial star_K\sigma
\]
and
\[
(K\setminus Intstar_K\sigma) \cap(cone\partial star_K\sigma)=\partial star_K\sigma,
\]
so our definition coincides with that in \cite{BP15}.

  Recall that $K \setminus Intstar_K \sigma = \{\tau \in K | \sigma \not\subset \tau \}$
 is a simplicial complex because of for every $\tau' \subset \tau$ and $\sigma \not\subset \tau$, one has $\sigma \not\subset \tau'$,
 then as $\tau \in K$ is a simplex of $K \setminus Intstar_K \sigma$ and $K$ is a simplicial complex,
 $\tau' \subset \tau$ is a simplex of $K \setminus Intstar_K \sigma$.

{\bf Theorem 2.9} {\it Let $\mathbb{P}=\{\sigma_1, \sigma_2, \cdots, \sigma_s\}$ be a simplicial
complement of $K$. Then $\{\mathbb{P}, \sigma, (\mathbb{P}-\sigma)\ast(m+1)\}$ is a simplicial complement
of $ss_{\sigma}K$ on the vertex set $[m+1]$, where
\begin{align*}
(\mathbb{P}-\sigma)\ast(m+1) = &
   \{(\sigma_1\setminus\sigma,m+1),(\sigma_2\setminus\sigma, m+1),\cdots,(\sigma_{s}\setminus\sigma,m+1)\}.
\end{align*} }

{\bf Proof:} First we prove that $\{\sigma,\mathbb{P}-\sigma\}$ is a simplicial complement of
$cone\partial star_K\sigma$ on the vertex set
$[m+1]=\{1,2,\cdots, m, m+1\}$.

From Lemma $2.6$ we know that $\mathbb{P} - \sigma$ is a simplicial complement of $star_{K}\sigma$
on the vertex set $[m]$.

A simplex $\tau$ on vertex set $[m]$ is not a simplex of
$\partial star_K\sigma=star_K\sigma\setminus Intstar_K\sigma$
if and only if $\tau \in Intstar_{K}\sigma$ or $\tau \notin star_{K}\sigma$, i.e.
$\sigma\subset \tau$ or there exists a $\sigma_i\setminus\sigma$ such that
$\sigma_i\setminus\sigma\subset \tau$, so
\[
2^{[m]}\setminus \partial star_K \sigma = 2^{[m]}\setminus K_{\{\sigma , \mathbb{P} - \sigma\}}([m]),
\]
$\{\sigma, \mathbb{P} - \sigma\}$
is a simplicial complement of
$\partial star_K\sigma$ on the vertex set $[m]$.

  Take the cone of $\partial star_K\sigma$ on the vertex
set $[m+1]$, a simplex $\tau\subset [m]$ or $(\tau, m+1)\subset [m+1]$ is not a simplex of
$cone \partial star_K\sigma=(m+1)\ast \partial star_K\sigma$ if and only
if $\tau$ is not a simplex of $\partial star_K\sigma$, i.e.
 $$2^{[m+1]}\setminus cone\partial star_K \sigma= 2^{[m+1]}\setminus K_{\{\sigma , \mathbb{P} - \sigma\}}([m+1]) ,$$
$\{\sigma, \mathbb{P} - \sigma\}$
is a simplicial complement of  $cone \partial star_K\sigma$ on the vertex set $[m+1]$.

  Second, we prove that $\{\mathbb{P}, \sigma, (m+1)\}$ is a simplicial complement of $K \setminus Intstar_K\sigma$ on the vertex set $[m+1]$.

  A simplex $\tau$ on the vertex set $[m]$ is not a simplex of $K\setminus Intstar_{K}\sigma$
if and only if  $\tau \notin K$ or $\tau \in Intstar_K\sigma$, i.e. there exists a $\sigma_i\in\mathbb{P}$ such that
$\sigma_i\subset \tau$ or $\sigma\subset\tau$.
$\{\mathbb{P}, \sigma\}$ is a simplicial complement of $K \setminus Intstar_K\sigma$ on the
vertex set $[m]$.

  Consider $K\setminus Intstar_K\sigma$ as a simplicial complex on the vertex set $[m+1]$,
$(m+1)$ does not appear in $K\setminus Intstar_K\sigma$.
It is a ghost vertex and $(m+1)$ is a missing face. So
\[
\{\mathbb{P}, \sigma, (m+1)\}
\]
is a simplicial complement of $K\setminus Intstar_K\sigma$ on the vertex set $[m+1]$.

  From Lemma $2.7$, we know that
$\{\mathbb{P}, \sigma, (m+1)\} \ast \{\sigma, \mathbb{P} - \sigma\}$ is a simplicial complement of
$ss_{\sigma}K=(K\setminus Intstar_K\sigma)\cup(cone\partial star_K\sigma)$, where
  \begin{align*}
\{\mathbb{P}, \sigma, (m+1)\} \ast \{\sigma, \mathbb{P} - \sigma\}
  = & \left\{\begin{array}{lllll}
     \mathbb{P} \ast \sigma, & \mathbb{P} \ast  \{\mathbb{P}- \sigma\}, \\
     \sigma \ast \sigma, & \sigma \ast \{\mathbb{P} - \sigma\}, \\
     (m+1) \ast \sigma,  & (m+1) \ast \{\mathbb{P} - \sigma \}
            \end{array}\right\}.
  \end{align*}

At last we complete the proof by showing that the simplicial complement
$\{\mathbb{P}, \sigma, (m+1)\} \ast \{\sigma, \mathbb{P} - \sigma\}$
is equivalent to $\{\mathbb{P}, \sigma, \{\mathbb{P}-\sigma\}\ast(m+1)\}$, i.e.
\begin{align*}
ss_{\sigma}K &= K_{\{\mathbb{P},\sigma,(m+1)\}\ast\{\sigma,\mathbb{P}-\sigma\}}([m+1])
 = K_{\{\mathbb{P}, \sigma, \{\mathbb{P}-\sigma\}\ast(m+1)\}}([m+1]),
\end{align*}

  First,
\[
\sigma\ast\sigma=\sigma\in \{\mathbb{P}, \sigma, (m+1)\} \ast \{\sigma, \mathbb{P} - \sigma\}.
\]
Every subset $\sigma_i \cup \sigma \in \mathbb{P} \ast \sigma$, $(\sigma,m+1)\in (m+1) \ast \sigma$
and $\sigma \cup (\sigma_i \setminus \sigma) \in \sigma\ast \{\mathbb{P}-\sigma\}$ contain $\sigma$.
They could be removed from $\{\mathbb{P}, \sigma, (m+1)\} \ast \{\sigma, \mathbb{P} - \sigma\}$, so
  \begin{align*}
\{\mathbb{P}, \sigma, (m+1)\} \ast \{\sigma, \mathbb{P} - \sigma\}
 \simeq  \left\{\begin{array}{lllll}
       \mathbb{P} \ast  \{\mathbb{P}- \sigma\}, \\
      \sigma,  \\
      (m+1) \ast \{\mathbb{P} - \sigma \}
            \end{array}\right\}.
  \end{align*}

  Then for any $\sigma_i \in \mathbb{P}$, one has $\sigma_i \setminus \sigma \in \mathbb{P}-\sigma$. So
\[
\sigma_i= \sigma_i \cup (\sigma_i \setminus \sigma) \in \mathbb{P} \ast \{\mathbb{P}-\sigma\}.
\]
Any other $\sigma_i \cup (\sigma_j \setminus\sigma) \in \mathbb{P} \ast \{\mathbb{P}-\sigma\}$ contains
$\sigma_i$, they could be removed from $\mathbb{P} \ast \{\mathbb{P}-\sigma\}$.
Thus $\mathbb{P} \ast \{\mathbb{P}-\sigma\}$ is equivalent to $\mathbb{P}$ and
$\{\mathbb{P}, \sigma, (m+1)\} \ast \{\sigma, \mathbb{P} - \sigma\}$ could be
reduced to
\begin{align*}
    & \{\mathbb{P}, \sigma, \{\mathbb{P}-\sigma\}\ast(m+1)\} \\
  = & \{\mathbb{P}, \sigma, (\sigma_1\setminus \sigma,m+1),(\sigma_2\setminus \sigma,m+1),\cdots,(\sigma_{s}\setminus \sigma,m+1)\}.
\end{align*}
The Theorem follows. \qb

{\bf Remark:} If $\sigma$ is  not a simplex of $K$, we still have $\{\mathbb{P}, \sigma, \{\mathbb{P}-\sigma\}\ast(m+1)\}$
as a simplicial complement of a simplicial complex $ss_\sigma K$. In that case, there exists a $\sigma_i\in\mathbb{P}$ such that
$\sigma_i\subseteq\sigma$. So $\sigma$ could be removed from $\{\mathbb{P}, \sigma, \{\mathbb{P}-\sigma\}\ast(m+1)\}$
and $\sigma_i\setminus \sigma=\emptyset \in \mathbb{P}-\sigma$. Thus
$(\sigma_i\setminus \sigma, m+1)=(m+1)\in \{\mathbb{P}-\sigma\}\ast(m+1)$ and all the other
$(\sigma_j\setminus \sigma, m+1)$ could be removed from  $\{\mathbb{P}, \sigma, \{\mathbb{P}-\sigma\}\ast(m+1)\}$.
That is to say that $(m+1)$ is a missing face and
\[
\{\mathbb{P}, \sigma, \{\mathbb{P}-\sigma\}\ast(m+1)\}\simeq \{\mathbb{P}, (m+1)\}
\]
is still a simplicial complement of $ss_\sigma K=K$ but on the vertex set $[m+1]$ and a ghost vertex $(m+1)$ is added.

We still call it the stellar subdivision at $\sigma$ on $K$. \qb

{\bf Example 3} In Example 1, we make stellar subdivision at $\sigma=(1,2)$ on $K$ (see Figure 2)

\begin{center}
\setlength{\unitlength}{0.7mm}
\begin{picture}(150,35)(-33,-15)
\put(0,20){\line(0,-1){40}}
\put(20,0){\line(-1,1){20.2}}
\put(20,0){\line(-1,-1){20.2}}
\put(-40,-10){\line(81,61){40}}
\put(-40,-10){\line(83,-21){40}}
\put(-40,-10){\line(6,1){25}}
\put(0,20){\line(-15,-26){15}}
\put(0,-20){\line(-20,19){15}}
\put(0,22){$1$}
\put(2,-22){$2$}
\put(21,0){$3$}
\put(-43,-10){$4$}
\put(-17,-11){$5$}

\put(35,0){$\stackrel{ss_{\sigma}K}{\Longrightarrow}$}
\put(90,20){\line(-2,-10){4}}
\put(90,-20){\line(-2,10){4}}
\put(110,0){\line(-1,1){20.2}}
\put(110,0){\line(-1,-1){20.2}}
\put(50,-10){\line(81,61){40}}
\put(50,-10){\line(83,-21){40}}
\put(50,-10){\line(6,1){25}}
\put(90,20){\line(-15,-26){15}}
\put(90,-20){\line(-20,19){15}}
\put(86,0){\line(-40,-21){11}}
\dashline{3}(50,-10)(86,0)
\put(90,22){$1$}
\put(92,-22){$2$}
\put(111,0){$3$}
\put(47,-10){$4$}
\put(73,-11){$5$}
\put(88,0){$6$}

\end{picture}
\end{center}

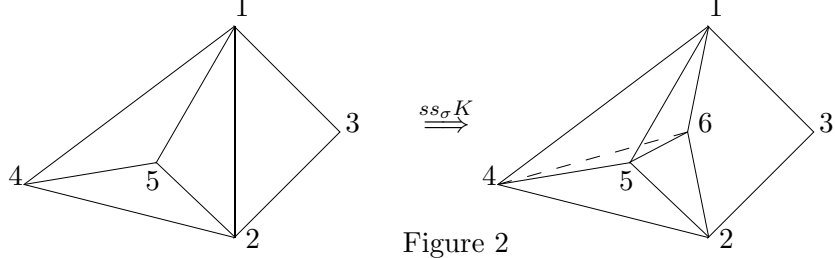
\centerline{Figure 2}

  $\mathbb{P}=\{(1,2,4,5), (1,2,3), (3,4), (3,5)\}$ is a simplicial complement of $K$, $\sigma=(1,2)$, so
\begin{align*}
\{\mathbb{P}-\sigma\}\ast(6)=\{(4,5), (3), (3,4), (3,5)\}\ast (6)
= & \{(4,5,6), (3,6), (3,4,6), (3,5,6)\}.
\end{align*}
\[
\{\mathbb{P}, \sigma, \{\mathbb{P}-\sigma\}\ast (6)\}
  =  \left\{\begin{array}{lllll}
  (1,2,4,5), (1,2,3), (3,4), (3,5),\\
   (1,2)=\sigma,\\
   (4,5,6), (3,6), (3,4,6), (3,5,6)
  \end{array}\right\}
\]
is a simplicial complement of $ss_\sigma K$.

\section{Construction}

  After given the simplicial complement of stellar subdivision, we construct our moment-angle manifolds
whose cohomology has torsion.

{\bf Lemma 3.1} {\it Let $K$ be a simplicial complex on the vertex set $[m]$ and
$$\mathbb{P}=\{\sigma_1, \sigma_2, \cdots, \sigma_s\}$$ be a simplicial complement of it.
Let $I$ be a subset of the vertex set $[m]$. Then
\[
\mathbb{P}|_I=\{\sigma_i \in \mathbb{P}| \sigma_i \subset I\}
\]
is a simplicial complement of the full subcomplex $K|_I$ on the vertex set I.}

{\bf Proof:} From its definition, we know that the full subcomplex
\[
K|_I=\{\sigma\in K |\sigma\subset I\}
\]
is a simplicial complex on the vertex set $I$.
A subset $\tau$ on the vertex set $I$ is not a simplex of $K|_I$ if and only if $\tau$
is not a simplex of $K$. i.e. there exists a non-face $\sigma_i \in \mathbb{P}$
such that $\tau_i \subset \tau$. Note that $\tau\subset I$, $\tau_i\subset\tau\subset I$.
The Lemma follows. \qb

{\bf Theorem 3.2 \space (Construction)}
 {\it Let $K$ be a subcomplex (not a full subcomplex) of a simplicial sphere $L_0$ on the vertex set $[m]$,
$\mathbb{M} = \{\sigma_1,\sigma_2,\cdots,\sigma_{s}\}$ be the set of missing faces of $K$, which are also simplices of $L_0$.
As following, make stellar subdivisions at $\sigma_1,\sigma_2,\cdots,\sigma_{s}$ on $L_0$ one by one
\[\begin{array}{cccccccccccc}
L_1 =ss_{\sigma_1}L_0,  & L_2 =ss_{\sigma_2}L_1, & \cdots,&  L_{s} =ss_{\sigma_{s}}L_{s-1}.
\end{array}\]
Then $K$ becomes a full subcomplex of $L_s$, $K = L_s|_{[m]}$.}

{\bf Proof:}
Let $\mathbb{P}_0=\{\tau_1, \tau_2, \cdots, \tau_r\}$ be a simplicial complement of $L_0$ on $[m]$.
From Theorem $2.9$ we know that
\[
\mathbb{P}_1=\{\mathbb{P}_0, \sigma_1, \mathbb{P}_1'\}
\]
is a simplicial complement of $L_1=ss_{\sigma_1}L_0$ on $[m+1]$,
where $$\mathbb{P}_1^{'}=\{\mathbb{P}_0-\sigma_1\}\ast (m+1).$$

  By induction, we get a simplicial complement of $L_s=ss_{\sigma_s}L_{s-1}$ on $[m+s]$ as
\begin{align*}
\mathbb{P}_s = & \{\mathbb{P}_{s-1}, \sigma_s, \mathbb{P}_{s}'\}\\
             = & \{\mathbb{P}_0, \sigma_1, \sigma_2, \cdots, \sigma_s,
                 \mathbb{P}_1', \mathbb{P}_2', \cdots, \mathbb{P}_s'\},
\end{align*}
where $$\mathbb{P}_{i}'=\{\mathbb{P}_{i-1}-\sigma_{i}\}\ast (m+i).$$

Note that every non-face in
$\mathbb{P}_i'$ contains $(m+i)$ as a vertex. From Lemma $3.1$ we know that
\[
\mathbb{P}_s|_{[m]} =\{\mathbb{P}_0, \sigma_1, \sigma_2, \cdots, \sigma_s\}
\]
is a simplicial complement of the full subcomplex $L_s|_{[m]}$.

Finally, we consider the simplicial complement $\mathbb{P}_s|_{[m]}$. Note that $K$
is a subcomplex of $L_0$, every non-face $\tau_i\in \mathbb{P}_0$ is not a simplex
of $K$, so there exists a $\sigma_j\in \mathbb{M}$ such that $\sigma_j\subseteq\tau_i$.
Then $\tau_i$ could be removed from $\{\mathbb{P}_0, \sigma_1, \sigma_2, \cdots, \sigma_s\}$.

Thus
\[
\mathbb{P}_s|_{[m]} =\{\mathbb{P}_0, \sigma_1, \sigma_2, \cdots, \sigma_s\}
 \simeq  \{\sigma_1, \sigma_2, \cdots, \sigma_s\}=\mathbb{M}
\]
which is the set of missing faces of $K$. The theorem follows. \qb

{\bf Remark:}
If $L_0$ is also a polytopal sphere, the stellar subdivision of $L_0$ is also polytopal.
It has been proved in a geometric sense by G.~Ewald and G.~C.~Shephard in
\cite{ES74}.

Let $\widetilde{L}_0$ be the simplicial polytope and its boundary
$\partial \widetilde{L}_0=L_0$ be the polytopal sphere.
If $\sigma$ is a simplex of $L_0$ and $\sigma$ is the intersection of the facets
(maximal simplices of $L_0$)
$F_{i_1}, F_{i_2}, \cdots , F_{i_r}$, one can take any point ${p}$ beyond the
facets $F_{i_1}, F_{i_2}, \cdots , F_{i_r}$ and beneath the other facets
 (See \cite{Gr03} p.78 for the definitions of ¡±beyond¡± and ¡±beneath¡±).
The stellar subdivision $ss_{\sigma}\partial \widetilde{L}_0$
is the boundary of the convex hall of $\widetilde{L}'_0=conv(\widetilde{L}_0 \cup {p})$.

It could also be proved from the duality of polytopes.

Let $\widetilde{L}_0$ be the simplicial polytope corresponding to $L_0$, and
$P_0$ be the dual simple polytope, (the vertex of $L_0$ corresponding to the facet
while the facet of $L_0$ corresponding to the vertex of $P_0$). Let $\sigma=(i_1, i_2, \cdots, i_k)$ be a simplex
of $L_0$, make a stellar subdivision at $\sigma$ on $L_0$
is equivalent, though the duality of polytopes,  to cutting off the face
$\sigma^*=F_{i_1}\cap F_{i_2}\cap \cdots \cap F_{i_k}$
in $P_0$ by a generic hyperplane.
The cutting off operation on a simple polytope is still simple, so $ss_{\sigma}\partial \widetilde{L}_0$ is polytopal. \qb

\section{Application}

{\bf Proposition 4.1}
{\it The cohomology of
differentiable moment-angle manifolds may have torsion of any order.}

{\bf Proof:} Let $L_0$ be a polytopal sphere and $K$ be a subcomplex of $L_0$, whose cohomology has torsion.
Construct a new polytopal sphere $L_s$ by Theorem $3.2$, then $K$ becomes a full subcomplex of $L_s$, while both $\mathbb{R}\mathcal{Z}_{L_s}$ and $\mathcal{Z}_{L_s}$ are framed
differentiable manifolds. From Hochster's Theorem, the cohomology of $\mathbb{R}\mathcal{Z}_{L_{s}}$ and $\mathcal{Z}_{L_{s}}$
has $\widetilde{H}^*(K)$ as a summand and then have torsion.

At least, every simplicial complex $K$ with $m$ vertexes is a subcomplex of the polytopal sphere $\partial\Delta^{m-1}$.
So the cohomology of differentiable moment-angle manifolds could have any torsion. \qb

Here is a example.

{\bf Example 4}
Let $K$ be the triangulated $mod$ $3$ Moore space (see Figure 3) which can be embedded in $6$-dimensional
polytopal sphere
\[
L_0=\partial \Delta^7=\partial(1,2,3,4,5,6,7,8).
\]

\begin{center}
\setlength{\unitlength}{0.7mm}
\begin{picture}(100,60)(-50,-25)
\put(0,0){\circle{50}}
\put(7.5,10){\line(-1,0){15}}
\put(-7.5,10){\line(-1,-2){6.5}}
\put(7.5,10){\line(1,-2){6.5}}
\put(-0.5,-13){\line(14,10){14.5}}
\put(-0.5,-13){\line(-13,10){13.4}}
\put(-0.5,-13){\line(8,23){8}}
\put(7.5,10){\line(-210,-127){21.5}}
\put(0,25){\line(-1,-2){10}}
\put(0,25){\line(1,-2){10}}
\put(16.10,19.15){\line(-12,-13){8.5}}
\put(24.62,4.35){\line(-191,63){17}}
\put(24.62,4.35){\line(-3,-2){11}}
\put(21.65,-12.5){\line(-11,14){8}}
\put(8.55,-23.5){\line(11,42){5.5}}
\put(8.55,-23.5){\line(-90,105){9}}
\put(-8.55,-23.5){\line(90,117){8.2}}
\put(-21.65,-12.5){\line(40,-1){21}}
\put(-21.65,-12.5){\line(33,41){8}}
\put(-24.62,4.35){\line(191,63){17}}
\put(-24.62,4.35){\line(3,-2){11}}
\put(-16.10,19.15){\line(12,-13){8.5}}
\put(22.65,-14.5){$1$}
\put(-1,26){$1$}
\put(20.15,18){$2$}
\put(25.5,3){$3$}
\put(8.5,-27.5){$2$}
\put(-10.5,-27.5){$3$}
\put(-27,-14.5){$1$}
\put(-28.5,3){$2$}
\put(-22.15,18){$3$}
\put(6.5,12){$4$}
\put(13,-2){$5$}
\put(-3,-11){$6$}
\put(-16.5,-2){$7$}
\put(-9.5,12){$8$}
\end{picture}
\end{center}

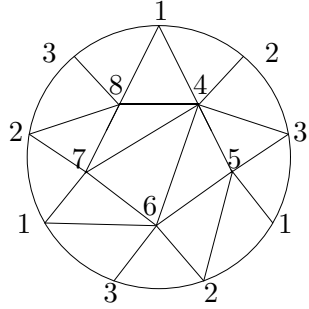
\centerline{Figure 3}

The set of missing faces of $L_0$ is
\[
\mathbb{P}_0=\{(1,2,3,4,5,6,7,8)\}.
\]
The set of missing faces of $K$ is
\[
\mathbb{M}=
\left\{\begin{array}{llllllllllllll}
(1,2,3),(1,2,6),(1,2,8),(1,3,4),\\
(1,4,5),(1,4,6),(1,4,7),(1,5,6),(1,7,8),\\
(2,3,5),(2,4,5),(2,4,6),(2,4,7),(2,4,8),(2,6,7),\\
(3,4,6),(3,4,8),(3,5,6),\\
(3,7), (5,8),(5,7), (6,8)\\
\end{array}\right\}\tag{4.2}
\]
and the set of maximal simplices of K is
\[
\left\{\begin{array}{llllllllllllll}
(1,2,4),(1,2,5),(1,2,7),(1,3,5),\\
(1,3,6),(1,3,8),(1,4,8),(1,6,7),\\
(2,3,4),(2,3,6),(2,3,8),(2,5,6),(2,7,8),\\
(3,4,5),(4,5,6),(4,6,7),(4,7,8)
\end{array}\right\}.
\]

Making $22$ stellar subdivisions at missing faces of $K$ on $\partial\Delta^7$, we
thus obtain a $6$-dimensional polytopal sphere $L_{22}$ with $30$ vertices which has $K$
as a full subcomplex. The real moment-angle manifold corresponding to $L_{22}$ is of $6$-dimensional
 while the complex one is of $37$-dimensional where $H^3(\mathbb{R}\mathcal{Z}_{L_{22}})$ and $H^{11}(\mathcal{Z}_{L_{22}})$ has $\widetilde{H}^2(K)=\mathbb{Z}/3$ as a summand.

Passing to the dual, $\Delta^7$ is the dual simple polytope of $\partial\Delta^7$
with facets numbered as vertexes of $\partial\Delta^7$. Making stellar subdivision on $\partial\Delta^7$ at $\sigma=(i_1,i_2,\cdots,i_r)$ is dual to
cutting off face $F_{i_1} \cap F_{i_2} \cap \cdots \cap F_{i_r}$ in $\Delta^7$,
$$
\xymatrix{
K \ar@{^{(}->}[r] & \partial \Delta^7 = \partial (\Delta^{7^{*}}) \ar@{<=>}[d] \ar[r]^-{\mbox{s.s.}} & L_{22} \ar@{<=>}[d]\\
& \Delta^7 \ar[r]^{\mbox{cut off}} & P_{22}.
}
$$
After cutting off the faces $F_{i_1}\cap F_{i_2}\cap \cdots \cap F_{i_r}$ numbered at $\mathbb{M}$ in (4.2),
one get a simple polytope $P_{22}$. The cohomology of the moment-angle manifold
corresponding to $P_{22}$ has $H^2(K)=\mathbb{Z}/3$ as a summand and then has torsion.
If we only cut off $\{1, 2,\cdots, 8\} \setminus \sigma$ for every maximal simplex $\sigma$ of $K$ in $\Delta^7$
as F.Bosio and L.Meersseman did in \cite{BM06} (Theorem $11.12$), we do not get torsion.

 Compute the missing faces after making stellar subdivision at $(1,2,3)$ and $(3,7)$ on $\partial\Delta^7$ in different sequence,
one has
\begin{enumerate}
\item We make stellar subdivision at $(1,2,3)$ on $L_0=\partial\Delta^7$ at first, then make stellar subdivision at $(3,7)$.

From Theorem $2.10$ we know that,
\begin{align*}
\mathbb{P}_0= & \{\underline{(1,2,3,4,5,6,7,8)}\}\\
\sigma_1= & (1,2,3)\\
(\mathbb{P}_0-\sigma_1)\ast (9) = & \left\{(4,5,6,7,8,9)\right\}
\end{align*}
is a simplicial complement of $L_1=ss_{(1,2,3)}L_0$. After removing the larger non-faces
$(1,2,3,4,5,6,7,8)$, we get the set of missing faces of $L_1$
\[
\mathbb{M}_1=\left\{(1,2,3),(4,5,6,7,8,9)\right\}.
\]

Then we make stellar subdivision at $(3,7)$ on $L_1$ and get the set of missing faces of $L_2=ss_{(3,7)}L_1$
\[
\mathbb{M}_2=\left\{(1,2,3),(4,5,6,7,8,9),(3,7),(1,2,10),(4,5,6,8,9,10)\right\}.
\]

\item Similarly we make stellar subdivision at $(3,7)$ on $L_0$ at first, then make stellar subdivision at $(1,2,3)$,
the resulting set of missing faces of $L_2 '$ is
\[
\mathbb{M}_2 '=\left\{(3,7),(1,2,4,5,6,8,9),(1,2,3),(7,10),(4,5,6,8,9,10)\right\}.
\]
\end{enumerate}

  It is easy to see that two simplicial complexes $K$ and $K'$ on vertex set $I$ are combinatorial
equivalent if and only if their sets of missing faces $\mathbb{M}$ and $\mathbb{M}'$ are equivalent,
i.e. there exists a one to one correspondence $\phi : I \rightarrow I$ that gives a one to one
 correspondence between $\mathbb{M}$ and $\mathbb{M}'$.

  Comparing with these two sequence, we can find that $L_2$ has one 2-vertex missing faces
$(3,7)$ while $L_2 '$ has two $(3,7),\ (7,10)$.
This implies that $L_2$ is not combinatorially isomorphic to $L_2 '$ and this difference might persist during the later
stellar subdivisions.

{\bf Remark:} Though $K$ will be a full subcomplex of $L_{s}$ in every sequence of making stellar subdivisions at $K$'s missing faces,
the combinatorial structure of $L_{s}$ may not be combinatorially isomorphic in different sequences. \qb

\end{document}